\documentclass[12pt]{amsart}
\begin{document}

\title{Liftez les Sylows!\\
Une suite \`a\\``Sous-groupes p\'eriodiques d'un groupe stable''}

\author{Bruno Poizat$^1$}
\address{\hskip-\parindent
        Bruno Poizat\\
        Institut Girard Desargues\\
        Universit\'e Claude Bernard (Lyon-1)\\
        b\^atiment 101 (math\'ematiques), 43, boulevard du 11 novembre 1918\\
        69621 Villeurbanne, France}
\email{poizat@desargues.univ-lyon1.fr}

\author{Frank O~Wagner$^2$}
\address{\hskip-\parindent
        Frank O~Wagner\\
        Mathematical Institute\\
        University of Oxford\\
        24--29 St Giles'\\
        Oxford OX1 3LB, UK}
\email{wagner@maths.ox.ac.uk}

\thanks{Research at MSRI is supported in part by NSF grant DMS-9701755.\\
{}$^1$ La unua geskribisto varmege dankas pro la \underline gia gastemo la Instituton
de Ser\underline cado de Matematikaj Sciencoj (MSRI) de la urbo Berkeley, en kiu
estis prilabora tiu artikoleto en la monato januaro de jaro mil na\underline u
cent na\underline u dek ok post la naski\underline go Krista.\\
{}$^2$ Heisenberg-Stipendiat der Deutschen Forschungsgemeinschaft (Wa 899/2-1).}

\begin{abstract} If $G$ is an omega-stable group with a normal definable subgroup
$H$, then the Sylow-$2$-subgroups of $G/H$ are the images of the
Sylow-$2$-subgroups of $G$.\\[4pt]
{\sc Zusammenfassung.} Sei $G$ eine omega-stabile Gruppe und $H$ ein
definierbarer Normalteiler von $G$. Dann sind die Sylow-$2$-Untergruppen von
G/H Bilder der Sylow-2-Untergruppen von $G$.\end{abstract}

\maketitle

Si $H$ est un sous-groupe normal d'un groupe fini $G$, les  $p$-sylows de
$G/H$  sont les images des  $p$-sylows de $G$ ; c'est une cons\'equence directe
de la conjugaison des sylows, et du fait que l'ordre d'un sylow est la
puissance maximale de $p$ divisant l'ordre du groupe ambiant.

Nous g\'en\'eralisons ici cette propri\'et\'e au cas o\`u $H$ est
d\'efinissable,  $p=2$, et $G$ est om\'ega-stable, ou m\^eme seulement
stable et menu, ou bien est stable et p\'eriodique; ou encore si $H$
est toujours d\'efinissable, $p$ premier quelconque, mais $G$ stable
et r\'esoluble par fini.

Ces hypoth\`eses nous ont permis dans \cite{pw93} de  montrer que
les sylows de $G$ (comme ceux de $G/H$!) sont localement finis et
conjugu\'es. Si nous renon\c cons \`a la seconde hypoth\`ese du Th\'eor\`eme
 14 (ii) de cet article, c'est que nous n'allons pas jusqu'\`a
pr\'etendre que  ${\mathbb Z}/2{\mathbb Z}$ se rel\`eve en un  $2$-sylow de $\mathbb Z$\!!

Nous faisons para\^\i tre cette note par remord de n'avoir pas fait cette
remarque alors (elle attire l'attention sur la difficult\'e de la
d\'emonstration du Th\'eor\`eme 14(ii) !), et aussi parce que ce r\'esultat
 simple fait d\'efaut dans un manuel sur le sujet comme \cite{bn94}.

Les cardinaux infinis manquant de valuation  $p$-adique, nous proc\'edons de
mani\`ere fort diff\'erente en commen\c cant par relever les  $p$-groupes
nilpotents. Nous montrons tout d'abord que si $G/H$ lui-m\^eme est un
$p$-groupe d'exposant fini, il se rel\`eve sur n'importe quel sylow  $S$
de $G$.

En effet, nous savons qu'alors $G/H$ est nilpotent; soit $H_1$ le
centre de $G$ modulo $H$: c'est \'egalement un groupe
d\'efinissable. Nos hypoth\`eses impliquent que tout \'el\'ement de
$G/H$ se rel\`eve en un $p$-\'el\'ement de $G$: il suffit de voir que
tout \'el\'ement $g$ de $G$ est congru modulo $H$ \`a un \'el\'ement
d'ordre fini; c'est \'evident dans le cas p\'eriodique; dans le cas
menu, cela vient de ce que $g$ est contenu dans un groupe ab\'elien
d\'efinissable $A$, et que $A\cap H$ se d\'ecompose en un groupe
divisible et un groupe d'exposant fini \cite[Lemme 13]{pw93}. Donc
tout \'el\'ement central dans $G/H$, se relevant dans un conjugu\'e de
$S$, a bien un ant\'ec\'edent dans $S$ lui-m\^eme. Cela r\`egle le cas
o\`u $G/H$ est commutatif; pour le cas g\'en\'eral, nous proc\'edons
par induction sur la classe de nilpotence de $G/H$, si bien que nous
savons que tout \'el\'ement de $G/H_1$ se rel\`eve dans $S$. Tout
\'el\'ement $g$ de $G$ s'\'ecrit donc comme produit d'un \'el\'ement de
$S$ et d'un \'el\'ement de $H_1$, ce dernier s'\'ecrivant comme
produit d'un \'el\'ement de $S$ et d'un \'el\'ement de $H$;
finalement, $g$ est bien congru modulo $H$ \`a un \'el\'ement de $S$.

Il nous suffit maintenant de montrer que tout sylow $\Sigma$ de $G/H$
s'exprime comme r\'eunion d'une suite croissante
$\Sigma_1,\ldots,\Sigma_n,\ldots$ de sous-groupes d\'efinissables
d'exposant fini; en effet, nous rel\`everons d'abord $\Sigma_1$ sur un
sylow $S_1$ de son image r\'eciproque, puis $\Sigma_2$ sur un sylow
$S_2$ de son image r\'eciproque contenant $S_1$, et ainsi de suite.

Changeant alors de notation, nous consid\'erons un sylow $S$ d'un
groupe $G$ satisfaisant nos hypoth\`eses. $S$ est localement fini et
nilpotent par fini; il est donc possible de trouver un sous-groupe
d\'efinissable nilpotent $G_1$ de $G$, normalis\'e par $S$, tel que
l'intersection $S_1$ de $S$ et de $G_1$ soit d'indice fini dans
$S$. Comme $S$ normalise l'unique sylow $S_2$ de $G_1$, il le
contient, et en fait $S_1=S_2$: $S_1$ est cet unique sylow de
$G_1$. D'apr\`es le Lemme 3 de \cite{pw93}, le quotient de $S_1$ par
son centre est d'exposant fini $p^a$. On trouve alors un entier $b$,
facilement calculable en fonction de $p^a$ et de la classe de
nilpotence de $S_1$, tel que pour tout $n\ge b$ les \'el\'ements de
$S_1$ d'ordre au plus $p^n$ forment un groupe $\Pi_n$ (voir
\cite[Lemma 1.2.20]{wa97}; dans le cas menu, on le voit de fa\c con
plus directe en observant que le centre de $S_1$ est engendr\'e par un
groupe d'exposant fini et un groupe divisible, si bien que $S_1$
lui-m\^eme est engendr\'e par un groupe d'exposant fini et un groupe
central divisible.) $\Pi_n$ est d\'efinissable car c'est en fait
l'ensemble des \'el\'ements d'ordre divisant $p^n$ de $G_1$. Notons
$\Gamma$ un groupe, fini, engendr\'e par un syst\`eme de
repr\'esentants de $S/S_1$; comme $\Pi_n$ est caract\'eristique dans
$S_1$, il est normalis\'e par $\Gamma$, et engendre avec ce dernier un
groupe d\'efinissable $\Sigma_n$, puisque $\Pi_n$ est d'indice fini
dans $\Sigma_n$. $S$ est bien l'union des $\Sigma_n$.

On observe que le r\'esultat de \cite{dw?} interdit aux $\Sigma_n$
d'\^etre tous normaux, et a fortiori caract\'eristiques, dans $S$
quand ce dernier n'est pas nilpotent.

\end{document}